\documentclass{article}

\usepackage{amsthm}
\usepackage{amssymb}
\usepackage{amsmath}
\usepackage{enumitem}
\usepackage[colorlinks=true,linkcolor=blue,citecolor=blue,urlcolor=blue]{hyperref}
\usepackage[all]{xy}
\usepackage{tikz-cd}

\def\tropfont#1{\mathsf{#1}}

\def\tropGm{\mathbf{G}_{\mathrm{trop}}}
\def\Gm{\mathbf G_m}
\def\logGm{\mathbf G_{\log}}
\def\Pic{\operatorname{Pic}}
\def\Hom{\operatorname{Hom}}

\def\Spec{\operatorname{Spec}}
\def\Sym{\operatorname{Sym}}
\def\LogPic{\operatorname{LogPic}}
\def\bLogPic{\operatorname{\mathbf{LogPic}}}
\def\bLogSch{\operatorname{\mathbf{LogSch}}}

\begin{document}

\title{Remarks on logarithmic \'etale sheafification}
\author{Sam Molcho and Jonathan Wise}
\date{\today}
\maketitle

\begin{abstract}
We prove criteria for a presheaf on logarithmic schemes to be a sheaf in the full logarithmic \'etale topology and describe several situations where the structure sheaf and logarithmic structure are logarithmic \'etale sheaves.  We deduce that the logarithmic Picard group is a stack in the full logarithmic \'etale topology on logarithmic schemes whose structure sheaves satisfy logarithmic \'etale descent.
\end{abstract}

\setcounter{secnumdepth}{5}
\renewcommand{\theparagraph}{\thesection.\arabic{paragraph}}

\section{Introduction}
The logarithmic \'etale topology is generated by strict \'etale covers, inherited from algebraic geometry, and \emph{logarithmic alterations}, which are compositions of logarithmic modifications and root stacks along Kummer extensions of the logarithmic structure.  Logarithmic alterations are universally surjective, proper, logarithmically \'etale monomorphims.  In conventional (that is, not logarithmic) algebraic geometry, surjective, proper, \'etale monomorphisms are necessarily isomorphims, but logarithmic alterations are far from being isomorphisms on the underlying schemes.  Nevertheless, many properties and invariants of logarithmic schemes are unchanged by replacement with a logarithmic alteration.  On the other hand, not everything is invariant under logarithmic alteration, and it can be a subtle matter to find a formulation that captures the invariance correctly.

The genesis of this paper is an error we committed when addressing this subtlety in a nearly final draft of \cite{logpic}.  In our analysis of logarithmic \'etale descent for logarithmic line bundles on a family of log curves $C \to S$, we incorrectly assumed invariance of the structure sheaf in the full logarithmic \'etale topology for an arbitrary logarithmic scheme $S$.  More precisely, we mistakenly assumed that, when $\tau : Y \to X$ is a logarithmic modification of logarithmic schemes, the induced homomorphism $\mathcal O_X \to \tau_\ast \mathcal O_Y$ is an isomorphism.  Examples constructed by K.\ Kato and C.\ Nakayama, summarized below, show that this homomorphism need not be injective or surjective.  In the published version of \cite{logpic}, we proved logarithmic \'etale descent for the logarithmic Picard group only under certain restrictive hypotheses on $S$.  In this note, we isolate this phenomenon by identifying the logarithmic schemes for which the above claim holds: they are exactly the logarithmic schemes for which the structure sheaf satisfies logarithmic \'etale descent. Our main result is that the original formulation of our descent statement can be salvaged if one works in the category of fine and saturated logarithmic schemes whose structure sheaves satisfy logarithmic \'etale descent:%
\footnote{We use the phrases `satisfies logarithmic \'etale descent' and 'is a sheaf in the logarithmic \'etale topology' synonymously.}

\paragraph{Theorem} \label{thm:main}
The fibered category of \'etale $\logGm$-torsors satisfies logarithmic \'etale descent on the category of fine and saturated logarithmic schemes whose structure sheaves satisfy logarithmic \'etale descent.

\paragraph{}
By \cite[Remark~4.4.3]{logpic}, logarithmic \'etale descent will always fail for logarithmic line bundles on logarithmic schemes whose structure sheaves do not satisfy logarithmic \'etale descent.

\paragraph{}
The demonstration of Theorem~\ref{thm:main} requires analysis of the sheafification of the structure sheaf of a logarithmic scheme in the logarithmic \'etale topology.  We establish criteria under which the structure sheaf satisfies logarithmic \'etale descent:

\paragraph{Theorem} A presheaf on logarithmically \'etale $X$-schemes is a sheaf in the full logarithmic \'etale topology if and only if it is a Kummer logarithmic \'etale sheaf and satisfies descent with respect to logarithmic modifications.  Moreover, the full logarithmic \'etale sheafification of a Kummer logarithmic \'etale sheaf $F$ is the universal presheaf receiving a morphism from $F$ and satisfying descent with respect to logarithmic modifications.

\paragraph{}
We furthermore single out common classes of logarithmic schemes whose structure sheaf satisfies log \'etale descent. For example: 

\paragraph{Theorem} Suppose $X$ is a logarithmic scheme that is logarithmically regular or which is logarithmically flat over a fine and saturated valuative%
\footnote{Fine logarithmic structures can be valuative if and only if they have rank~$1$.  However, this result should generalize to schemes with valuative logarithmic structures that have local charts but are not necessarily locally finitely generated.  We have not formulated it in that generality because we did not want to develop the background material necessary to work with such objects.}
	logarithmic scheme $S$. Then the structure sheaf of $X$ satisfies log \'etale descent. 

\paragraph{}
Apart from this criterion, we were only able to accomplish the most basic analysis of the logarithmic \'etale sheafification of the structure sheaf.  Many questions remain that we were unable to answer:

\subparagraph{Question}
Suppose that $X$ is a fine and saturated logarithmic scheme.  Let $X'$ be the spectrum of the logarithmic \'etale sheafification of the structure sheaf of $X$.  Is $X'$ of finite type over $X$?

\subparagraph{Question}
Can the logarithmic schemes whose structure sheaves are logarithmic \'etale sheaves be characterized in a more direct way?

\subparagraph{Question}
We know of several examples of fine and saturated logarithmic schemes whose structure sheaves do not satisfy logarithmic \'etale descent (see Section~\ref{sec:examples}).  What are the logarithmic \'etale sheafification of their structure sheaves?  We do not know how to calculate the logarithmic \'etale sheafification of the structure sheaf except in trivial examples.

\paragraph{}
Nevertheless, we believe that the condition is central in logarithmic geometry: it enters into consideration any time a problem is compactified by logarithmic techniques. 

\paragraph{Conventions}
If $S$ is a scheme with a logarthmic structure, we write $\varepsilon : M_X \to \mathcal O_X$ for the structural monoid homomorphism.  The monoid $M_X$ is written multiplicatively and the quotient $\bar M_X = M_X / \varepsilon^{-1}(\mathcal O_X^\ast)$ is written additively.  We write $M_X^{\rm gp}$ and $\bar M_X^{\rm gp}$ for the associated groups.  We write $\logGm$ and $\tropGm$, respectively, for the presheaves on logarithmic schemes whose values on $S$ are $\Gamma(S, M_X^{\rm gp})$ and $\Gamma(S, \bar M_X^{\rm gp})$.

If $X$ is a fine and saturated logarithmic scheme, a Kummer extension of the characteristic monoid of $X$ is an \'etale-locally finitely generated extension $\bar M_X \subset P$ such that, for every local section $\alpha$ of $P$, there is a positive integer $n$ such that $n \alpha \in \bar M_X$.  For any such extension, there is a universal algebraic stack $Y$ with integral, saturated logarithmic structure and projection $f : Y \to X$ such that $f^{-1} \bar M_X \to \bar M_Y$ extends to $f^{-1} P \to \bar M_Y$.  The logarithmic structure of $Y$ admits charts on a smooth cover.  The stack $Y$ is called the \emph{root stack} of $X$ along the extension $P$.  If, for every local section $\alpha$ of $P$, the integer $n$ can be chosen invertible over the domain of definition of $\alpha$ then $Y$ is called a root stack of \emph{invertible order}.

\paragraph{Acknowledgements}
We thank L.\ Herr and C.\ Nakayama for valuable discussions, as well as the anonymous referee who brought the error in \cite{logpic} to our attention. S.M was supported by the grant ERC-2017-AdG-786580-MACI. J.W.\ was supported by a Simons Foundation Fellowship, award \#822534, and Simons Travel Support, award \#636210.

\section{Saturated elements of the characteristic monoid}

\paragraph{Definition}
Let $\bar M$ be a fine, saturated, sharp monoid.  An element $\alpha \in \bar M^{\rm gp}$ will be called \emph{saturated} if $\alpha$ generates the associated group of every rank~$1$ localization of $\bar M$.

\paragraph{}

The idea is that if $X$ is a fine and saturated logarithmic scheme and $\alpha \in \Gamma(X, \bar M_X^{\rm gp})$ then $\alpha$ determines a morphism $X \to \tropGm$, and we should think of this as a saturated morphism if $\alpha$ satisfies the condition of the definition.

\paragraph{Lemma} \label{lem:saturated}
Let $\bar M$ be a fine and saturated monoid.  Suppose $\alpha \in \bar M^{\rm gp}$ is saturated.  Then $\bar M + \mathbf N \alpha \subset \bar M^{\rm gp}$ is saturated.

\begin{proof}
	See \cite[Lemma~4.4.11.1]{logpic}.
\end{proof}

\paragraph{Corollary} \label{cor:saturated}
Suppose $X$ is a fine and saturated logarithmic scheme and $\alpha$ is a saturated section of $\bar M_X^{\rm gp}$.  Let $Y$ be the universal \emph{fine} logarithmic scheme over $X$ where $\alpha \in \bar M_Y$.  Then $Y$ is saturated and its geometric fibers over $X$ are reduced.

\begin{proof}
It is immediate from Lemma~\ref{lem:saturated} that $Y$ is saturated.  We move on to showing that the geometric fibers are reduced.

	Writing $X[\alpha]$ for the universal fine logarithmic scheme $Y$ over $X$ such that $\alpha \in \bar M_Y$, if $x$ is a geometric point of $X$ with the induced logarithmic structure then $x \mathop\times_X X[\alpha] = x[\alpha]$.  It therefore suffices to demonstrate the corollary in the case where $X = x$.

We therefore assume that $X$ is the spectrum of an algebraically closed field.  There are now $3$ possibilities:  if $\alpha \in \bar M_X$ then $Y = X$ and we are done; if $-\alpha \in \bar M_X$ and $\alpha \not\in \bar M_X$ then $Y = \varnothing$ and we are done.

The final possibility is that neither $\alpha$ nor $-\alpha$ is in $\bar M_X$.  We claim that in this case the underlying scheme of $Y$ can be identified with the total space $Z = \Spec \Sym \mathcal O_X(-\alpha)$ of $\mathcal O_X(\alpha)$.  Indeed, the inequality $\alpha \geq 0$ on $Y$ corresponds to a section of $\mathcal O_Y(\alpha)$ and therefore a map $Y \to Z$.  On the other hand, let $\pi : Z \to X$ be the projection and define $M_Z$ to be the submonoid of $\mathcal O_Z^\ast \otimes_{\pi^{-1} \mathcal O_X^\ast} \pi^{-1} \bar M_X^{\rm gp}$ lying above $\pi^{-1} (\bar M_X + \mathbf N \alpha)$.  Let $\sigma : \mathcal O_Z(-\alpha) \to \mathcal O_Z$ be the tautological map.  Define $\varepsilon(a) = \sigma^n(a)$ when $a \in \mathcal O_Z^\ast(-n \alpha) \subset \bar M_Z$ and $\varepsilon(a) = 0$ for all other $a \in \bar M_Z$.  Since $\mathbf N \alpha$ is an ideal in $\bar M_X + \mathbf N \alpha$,%
\footnote{If $\beta \in \bar M_X + \mathbf N \alpha$ then $\beta = \gamma + n \alpha$ for some $\gamma \in \bar M_X$ and integer $n \geq 0$.  While this representation is not unique, if $\beta = m \alpha$ then we must have $\gamma = 0$.  Indeed, if $\beta = m \alpha$ then $\gamma = (m-n)\alpha$.  Since $\bar M_X$ is saturated and neither $\alpha$ nor $-\alpha$ lies in $\bar M_X$, we must therefore have $m = n$ and hence $\gamma = 0$.  It follows immediately from this observation that $\mathbf N \alpha$ is an ideal in $\bar M_X + \mathbf N \alpha$.}
	this is a monoid homomorphism, and one quickly verifies that it is a fine logarithmic structure.  Therefore the universal property of $Y$ gives a map $Z \to Y$.  Since $Y \to X$ and $Z \to X$ are monomorphisms, it follows that the maps $Z \to Y$ and $Y \to Z$ constructed above are inverse isomorphisms.  On the other hand $Z$ is reduced and its logarithmic struture is saturated by Lemma~\ref{lem:saturated}, so we may conclude.
\end{proof}

\paragraph{Lemma} \label{lem:root-sat}
Suppose that $X$ is a quasicompact, fine and saturated logarithmic scheme and $\alpha \in \Gamma(X, \bar M_X^{\rm gp})$.  Then there is a root stack $X'$ of $X$ and an integer $n > 0$ such that $n^{-1} \alpha$ is saturated on $X'$.

\begin{proof}
	See \cite[Lemma~4.4.11.5]{logpic}.
\end{proof}

\section{Fibers of toric morphisms}

\paragraph{Lemma}
Let $U$ be an affine toric variety with dense torus $j : T \to U$, equipped with its toric logarithmic structure.  The following groups are all isomorphic via the natural maps: 
\begin{enumerate}[label=(\arabic{*})]
\item the group of isomorphism classes of invertible sheaves on $[U/T]$;
\item the group of invertible sheaves on $[U/T]$, trivialized over the open point;
\item the group of invertible sheaves on $U$ trivialized over $T$;
\item $\Gamma(U, j_\ast (\mathcal O_T^\ast) / \mathcal O_U^\ast)$;
\item $\Gamma(U, \bar M_U^{\rm gp})$;
\item $\Gamma([U/T], \bar M_{[U/T]}^{\rm gp})$.
\end{enumerate}

\begin{proof}
We write $P_i$ for the $i$-th group in the list.  We have a map $P_2 \to P_1$ by forgetting the trivialization over the open point.  Since the fiber over the open point is always trivial, this is surjective.  Since all trivializations of an invertible sheaf $L$ over the open point are related by global isomorphisms of $L$, it is also injective.  Thus $P_1 \simeq P_2$.

By definition of the logarithmic structure on a toric variety $M_U^{\rm gp} = j_\ast \mathcal O_T^\ast$, so $P_4 \simeq P_5$.

If $L$ is an $\mathcal O_U^\ast$-torsor on $U$ with a trivialization $j^\ast L \simeq \mathcal O_T^\ast$ then we get an inclusion $L \to j_\ast \mathcal O_T^\ast$.  Thus $L$ is a $\mathcal O_U^\ast$-coset in $j_\ast \mathcal O_T^\ast$, hence corresponds to a section of $P_4$.  This shows $P_3 \simeq P_4$.

Note that $\bar M_U^{\rm gp}$ is pulled back from $\bar M_{[U/T]}^{\rm gp}$.  This induces an isomorphism $P_6 \simeq P_5$ because $\bar M_{[U/T]}^{\rm gp}$ is an \'etale sheaf and $U \to [U/T]$ is a universal submersion (it is surjective and flat).

	Finally, we have map $P_2 \to P_3$ by pullback.  For any flat scheme over $f : V \to [U/T]$, define $P(V)$ to be the group of invertible sheaves on $V$ trivialized on the preimage of the open point.  Then $P([U/T]) = P_2$ and $P(U) = P_3$.  It is immediate that $P$ is a sheaf.  But $U \mathop\times_{[U/T]} U \simeq U \times T$ and $P(U \times T) \simeq P(U)$ (by the isomorphism $P_3 \simeq P_6$).  Both the projection and the action map induce the same isomorphism $P(U) \to P(U \times T)$ (one can also see this from the isomorphism $P_3 \simeq P_6$) so by flat descent, $P_2 \simeq P_3$.
\end{proof}

\paragraph{Lemma} \label{lem:compatible-characters}
Suppose that $U$ is an affine toric variety.  Let $T \subset U$ be the dense torus.  Suppose that $V$ and $W$ are torus orbits in $U$, with $V$ in the closure of $W$.  If $L$ is an equivariant invertible sheaf on $U$ then the restrictions of $L$ to $V$ and to $W$ correspond to characters $\chi_V$ and $\chi_W$ in the character groups $X_V$ and $X_W$ of the stabilizers of $V$ and $W$.  The image of $\chi_V$ under the homomorphism $X_V \to X_W$ is $\chi_W$.

\begin{proof}
This is immediate: the map $X_V \to X_W$ is induced from the inclusion of the stabilizer subgroups $T_W \subset T_V$.
\end{proof}

\paragraph{Lemma} \label{lem:pic-splitting}
Let $X$ be a toric variety with dense torus $T$.  Let $\tropfont X = [X/T]$ and let $\tropfont Y$ be a closed substack.  The natural morphism $\Gamma(\tropfont Y, \bar M_{\tropfont Y}^{\rm gp}) \to \Pic(\tropfont Y)$ is a split injection.

\begin{proof}
An invertible sheaf on $\tropfont Y$ gives a section of $\bar M_{\tropfont Y}^{\rm gp}$ on each stratum, and these sections must be compatible by Lemma~\ref{lem:compatible-characters}.
\end{proof}

\paragraph{Lemma} \label{lem:characteristic-injective}
Let $X$ be an affine toric variety.  Let $\tau : Y \to X$ be a toric modification.  Let $x$ be a closed, torus invariant point of $X$.  Then $\Pic(Y) \to \Pic(\tau^{-1} x)$ is injective.

\begin{proof}
Let $\tropfont X$ and $\tropfont Y$ be the quotients of $X$ and $Y$ by their common dense torus, $T$.  Let $\tropfont V = [x/T]$ and let $\tropfont W = [ \tau^{-1} x / T ]$.  By Lemma~\ref{lem:pic-splitting}, the following homomorphism is a split injection:
\begin{equation*}
	\Pic(\tropfont Y) = \Gamma(\tropfont Y, \bar M_{\tropfont Y}^{\rm gp}) \xrightarrow{\sim} \Gamma(\tropfont W, \bar M_{\tropfont Y}^{\rm gp}) \to \Pic(\tropfont W)
\end{equation*}
Now use the following commutative diagram (with $W = \tau^{-1} x$):
\begin{equation} \label{eqn:a1} \vcenter{ \xymatrix{
\Pic(\tropfont Y) \ar[r] \ar[d] & \Pic(\tropfont W) \ar[d] \\
\Pic(Y) \ar[r] & \Pic(W)
} } \end{equation}
	The vertical arrow on the left is surjective, and its kernel is the character group of $T$.  Similarly, the $T$-linearizations of an invertible sheaf on $W$ form a torsor under $H^1( \mathrm B T, G )$ where $G$ is the algebraic group $G(Z) = \Gamma(Z \times W, \Gm)$.  Now, $W$ is proper and connected so $G \simeq \Gm \times U$ for some unipotent group $U$.  Since $U$ is an iterated extension of additive groups, $H^1( \mathrm B T, U ) = 0$, and therefore $H^1(\mathrm BT, G) = H^1(\mathrm BT, \Gm)$ is the character group of $T$.  Therefore the kernels of the vertical arrows in~\eqref{eqn:a1} are the same.  Since $\Pic(\tropfont Y) \to \Pic(Y)$ is surjective, it follows that $\Pic(Y) \to \Pic(W)$ is injective (by the snake lemma).
\end{proof}

\section{Logarithmic \'etale sheafification}

\paragraph{Lemma} \label{lem:log-etale-sheaf}
A presheaf on fine and saturated logarithmic schemes is a sheaf in the Kummer logarithmic \'etale topology if and only if it satisfies descent with respect to \'etale covers and root stacks of invertible order.%
\footnote{Counterintuitively, all root stacks satisfy the formal criterion for logarithmically \'etale morphsims, even though roots that are not coprime to the characteristic would appear to violate Kato's criterion for logarithmically \'etale morphisms \cite[Theorem~(3.5)]{Kato88}.  They do not violate the latter criterion because they are only logarithmic algebraic stacks and not logarithmic schemes.  

Given $\alpha \in \Gamma(X, \bar M_X)$, the $n$-th root of $X$ along $\alpha$ is the universal logarithmic algebraic stack $X'$ over $X$ over which $n^{-1} \alpha \in \Gamma(X', \bar M_{X'})$.  If $S \to S'$ is a strict infinitesimal extension of logarithmic schemes, then every commutative square
\begin{equation*} \begin{tikzcd}[ampersand replacement=\&]
	S \ar[r] \ar[d] \& X' \ar[d] \\
	S' \ar[r] \ar[ur,dashed] \& X
\end{tikzcd} \end{equation*}
has a unique lift.  Indeed, we have $\Gamma(S, \bar M_S) = \Gamma(S', \bar M_{S'})$, so a section of $\bar M_S$ has an $n$-th root if and only if the corresponding section of $\bar M_{S'}$ does.}
	It is a sheaf in the full logarithmic \'etale topology if, in addition, it satisfies descent with respect to any of the following classes of covers:
\begin{enumerate}[label=(\roman{*})]
\item logarithmic modifications;
\item logarithmic blowups;
\item pullbacks of $\mathbf P^1 \to \logGm$; or
\item pullbacks of $\mathbf P^1 \to \logGm$ along saturated maps.
\end{enumerate}

\begin{proof}
Since \'etale covers and root stacks of invertible order are Kummer logarithmic \'etale covers, and logarithmic modifications are logarithmic \'etale covers, it is clear that descent with respect to these classes is necessary to form a sheaf.  Conversely, suppose that $F$ is a presheaf that satisfies descent with respect to \'etale covers and root stacks of invertible order and let $R$ be a Kummer logarithmic \'etale covering sieve of a logarithmic scheme $X$.  We wish to show that any descent datum for $F$ relative to $R$ descends to $X$.  Since $F$ is a sheaf in the \'etale topology, the conclusion is an \'etale-local property of $X$, and we can therefore assume that $X$ has a global chart.

The sieve $R$ is generated by a collection of Kummer logarithmically \'etale morphisms $\{ Y_i \to X \}$, and we can assume that each $Y_i$ is quasicompact.  Since $X$ has a global chart and each $Y_i$ is quasicompact, each of the morphisms $Y_i \to X$ is the composition of an \'etale morphism (which is open) and a root stack (which is bijective on points and finite, hence closed, hence a homeomorphism), the image of each $Y_i \to X$ is open.  The images of the maps $Y_i \to X$ therefore form an open cover of $X$.  Since $X$ is quasicompact, we may assume that there are only finitely many $Y_i$ (replacing $R$ by a finer cover if necessary, which is harmless).  Each $Y_i \to X$ factors as $Y_i \to Y'_i \to X$ where $Y'_i$ is a root stack of invertible order of $X$.  The fine and saturated fiber product $Z' = Y'_1 \mathop\times_X \cdots \mathop\times_X Y'_n$ is also a root stack of invertible order, and the $Z_i = Y_i \mathop\times_{Y'_i} Z'$ form an \'etale cover of $Z'$.  A descent datum for $F$ relative to $R$ will descend to $Z'$ because $F$ is a sheaf in the \'etale topology, and then will descend to $X$ because $F$ satisfies descent with respect to root stacks of invertible order.

	The full logarithmic \'etale topology is generated by logarithmic blowups and Kummer logarithmic \'etale covers \cite[Proposition~3.9]{Nakayama17}.  Since logarithmic blowups are logarithmic modifications, logarithmic modifications certainly suffice as well.  Any logarithmic modification of $X$ can be refined, \'etale-locally on $X$, by a subdivision along hyperplanes, so subdivisions along hyperplanes also suffice.  Finally, if $\alpha \in \Gamma(X, \bar M_X^{\rm gp})$, there is a root stack $X'$ of $X$ and an integer $n > 0$ such that $n^{-1} \alpha$ lies in $\bar M_{X'}^{\rm gp}$ and $n^{-1} \alpha$ is saturated (Lemma~\ref{lem:root-sat}).
\end{proof}

\paragraph{Lemma} \label{lem:sheafification}
Let $X$ be a quasicompact, quasiseparated logarithmic scheme with a global chart.  Let $F$ be a sheaf in the Kummer logarithmic \'etale topology on the category of logarithmically \'etale $X$-schemes.  Write $F_Z$ for the restriction of $F$ to the small \'etale site of $Z$.  For each logarithmically \'etale $Y \to X$, let $G_Y = \varinjlim_{\rho} \rho_\ast F_Z$, with the colimit taken over all logarithmic modifications $\rho : Z \to Y$.  Then $G$ is the logarithmic \'etale sheafification of $F$.

\begin{proof}
Suppose that $H$ is a logarithmic \'etale sheaf on $X$ and $\rho : Z \to Y$ is a logarithmic modification.  Then $\rho$ is a logarithmically \'etale monomorphism, so $H(Y) \to H(Z)$ is an isomorphism.  If $F \to H$ is any morphism then $F(Y) \to H(Y)$ can be factored uniquely:
\begin{equation*}
F(Y) \to \rho_\ast F_Z(Y) = F(Z) \to H(Z) \xleftarrow{\sim} H(Y)
\end{equation*}
Thus $F(Y) \to H(Y)$ factors uniquely through $G(Y)$, so $F \to H$ factors uniquely through $G$.  It will therefore be sufficient to show that $G$ is a sheaf in the logarithmic \'etale topology.  According to Lemma~\ref{lem:log-etale-sheaf}, we need to check descent with respect to logarithmic modifications, root stacks of invertible order, and \'etale covers.

First, suppose that $\tau : Y \to X$ is a logarithmic modification.  Every logarithmic modification of $Y$ is a logarithmic modification of $X$ and these are cofinal in the filtered system of logarithmic modifications of $X$.  Therefore:
\begin{equation*}
	G_X = \varinjlim_{\rho : Z \to Y} \tau_\ast \rho_\ast F_Z = \tau_\ast \varinjlim_{\rho : Z \to Y} \rho_\ast F_Z = \tau_\ast G_Y
\end{equation*}
We are able to commute $\tau_\ast$ with the colimit because the colimit is filtered and $\tau$ is quasicompact and quasiseparated.

Now suppose that $\tau : Y \to X$ is a root stack.  Pullback along $\tau$ induces an equivalence between the filtered systems of logarithmic modifications of $Y$ and of $X$.  Therefore:
\begin{multline*}
	\tau_\ast G_Y = \tau_\ast \varinjlim_{\rho : Z \to X} \rho_\ast F_{Z \mathop\times_X Y} = \varinjlim_{\rho : Z \to X} \tau_\ast \rho_\ast F_{Z \mathop\times_X Y} \\
	= \varinjlim_{\rho : Z \to X} \rho_\ast \tau_\ast F_{Z \mathop\times_X Y} = \varinjlim_{\rho : Z \to X} \rho_\ast F_Z = G_X
\end{multline*}
We have used the quasicompactness and quasiseparatedness of $\tau$ to pass $\tau_\ast$ across the colimit and we have used the hypothesis that $F$ is a sheaf with respect to the Kummer logarithmic \'etale topology to get $\tau_\ast F_{Z \mathop\times_X Y} = F_Z$.

Finally, we must show that $G$ is an \'etale sheaf.  Suppose that $\{ Y_i \to X \}$ is an \'etale cover.  Refining the cover if necessary, we assume that there are finitely many $Y_i$ (by the quasicompactness of $X$) and that all $Y_i$ are quasicompact and quasiseparated.  Since $X$ has a global chart, every logarithmic modification of each $Y_i$ can be refined by the pullback of a logarithmic modification of $X$.  Since there are only finitely many $Y_i$, if $Z_i \to Y_i$ are logarithmic modifications of each of the $Y_i$ then there is a logarithmic modification $W \to X$ such that $W \mathop\times_X Y_i$ refines all of the $W_i$.  That is, the filtered system of `logarithmic modifications of the $Y_i$' is refined by the filtered system of `pullbacks of logarithmic modifications of $X$'.  If $\tau : R \to X$ is the sieve generated by the $Y_i$, we therefore have:
\begin{equation*}
	\tau_\ast G_R = \tau_\ast \varinjlim_{\rho : Z \to X} \rho_\ast F_{R \mathop\times_X Z} = \varinjlim_{\rho : Z \to X} \rho_\ast \tau_\ast F_{R \mathop\times_X Z} = \varinjlim_{\rho : Z \to X} \rho_\ast F_Z = G_X
\end{equation*}
Here we have used the quasisepartedness of $X$ to guarantee that $\tau$ is quasicompact and quasiseparated, so that we can commute the colimit with $\tau_\ast$, and we have used the assumption that $F$ is a sheaf in the \'etale topology to get $\tau_\ast F_{R \mathop\times_X Z} = F_Z$.
\end{proof}

We recall that the valuativization of a fine and saturated logarithmic scheme $X$ is the universal valuative logarithmic locally ringed space $X^{\rm val}$ equipped with a projection to $X$ \cite[Theorem~1.3.1]{Kato89}, \cite[3.12]{Nakayama17}.  Explicitly, $X^{\rm val} = \varprojlim_{Y \to X} Y$, with the limit taken over all logarithmic modifications $Y$ of $X$.

Suppose that $F$ is a presheaf on logarithmically \'etale logarithmic schemes over $X$ that is a sheaf in the Zariski topology.  Then $F$ induces a sheaf on $X^{\rm val}$:
\begin{equation*}
	F^{\rm val} = \varinjlim_{\varphi^{-1}} F_Y;
\end{equation*}
the colimit is taken over all logarithmic modifications $Y \to X$, with $\varphi : X^{\rm val} \to Y$ denoting the canonical projection.

C.\ Nakayama suggested the following reformulation of Lemma~\ref{lem:sheafification}:
\paragraph{Corollary}
	With the same hypotheses as in Lemma~\ref{lem:sheafification}, let $\pi : X^{\rm val} \to X$ denote the projection.  Then $\pi_\ast F^{\rm val}$ is the restriction of the sheafification of $F$ in the logarithmic \'etale topology to the Zariski site of $X$.%
	\footnote{The restriction to the Zariski site is included only to avoid defining the Kummer \'etale site of $X^{\rm val}$.}

\begin{proof}
	We observe first that, for any sheaf $G$ in the Zariski topology of $X$, the map $G \to \pi_\ast \pi^{-1} G$ is an isomorphism.  Indeed, it is sufficient to verify that the map of stalks $G_x \to \pi_\ast (\pi^{-1} G)_x$ is an isomorphism for each point $x$ of $X$.  But $\pi$ is surjective, closed, and quasicompact \cite[Proposition~3.13]{Nakayama17}.  It also has connected fibers, since all logarithmic modifications have connected fibers.  Therefore $\pi_\ast ( \pi^{-1} G)_x = \pi_\ast( \pi^{-1} G_x)$, where we abusively use $\pi$ for the restriction of $\pi$ to a fiber.

	Pushforward along $\pi$ also commutes with filtered colimits.  We conclude that
	\begin{equation*}
		\pi_\ast F^{\rm val} = \pi_\ast  \varinjlim \varphi^{-1} F_Y = \varinjlim \pi_\ast \varphi^{-1} F_Y = \varinjlim \tau_\ast F_Y ;
	\end{equation*}
	the colimits are taken over logarithmic modifications $\tau : Y \to X$. Thus $\pi_\ast F^{\rm val}$ is the logarithmic sheafification of $F$, by Lemma~\ref{lem:sheafification}.
\end{proof}

\section{Sheafification of the structure sheaf}

\paragraph{Corollary} \label{cor:logification}
Suppose that $X$ is a quasicompact, quasiseparated, fine, and saturated logarithmic scheme that has a global chart.  Then the logarithmic \'etale sheafifications of $\mathcal O_X$, $\mathcal O_X^\ast$, $M_X$, and $M_X^{\rm gp}$ are $\varinjlim \tau_\ast \mathcal O_Y$, $\varinjlim \mathcal O_Y^\ast$, $\varinjlim \tau_\ast M_Y$, and $\varinjlim \tau_\ast M_Y^{\rm gp}$, respectively, with the colimit taken over all logarithmic modifications $\tau : Y \to X$.

\begin{proof}
All of $\mathcal O_X$, $\mathcal O_X^\ast$, $M_X$, and $M_X^{\rm gp}$ are Kummer logarithmic \'etale sheaves~\cite[Theorems~3.1 and~3.2]{Kato21}, since $\mathcal O_X$ is representable by $\mathbf A^1$ with the trivial logarithmic structure, $\mathcal O_X^\ast$ is representable by $\Gm$ with the trivial logarithmic structure, and $M_X$ is representable by $\mathbf A^1$ with the toric logarithmic structure.  We conclude by Lemma~\ref{lem:sheafification}.
\end{proof}

\paragraph{Lemma} \label{lem:log-mod-pf}
Let $X$ be a (fine and saturated) logarithmic scheme and let $\tau : Y \to X$ be a logarithmic modification.  Set $\mathcal O'_X = \tau_\ast \mathcal O_Y$ and $M'_X = \tau_\ast M_Y$.  Let $X' = \Spec_X \mathcal O'_X$ with the logarithmic structure induced from $M'_X$.  Then we have the following properties:
\begin{enumerate}[label=(\roman{*})]
\item \label{item:7} $M'_X \to \mathcal O'_X$ is a logarithmic structure;
\item \label{item:8} the projection $\mu : X' \to X$ is finite;
\item \label{item:9} $\mu$ is a universal bijection (in the category of schemes);
\item \label{item:15} the residue field extensions of $\mu : X' \to X$ are trivial;
\item \label{item:10} $\mu^\ast$ is an equivalence between the strict \'etale sites of $X$ and $X'$;
\item \label{item:11} $\bar M_X \to \bar M'_X$ is an isomorphism;
\item \label{item:12} $\mu^\ast$ is an equivalence between the logarithmic \'etale sites of $X$ and $X'$;
\item \label{item:13} $\mathrm R\mu_\ast \mathcal O_{X'} = \mathcal O'_X$ and $\mu_\ast M_{X'} = M'_X$ and $\mathrm R \mu_\ast M_{X'}^{\rm gp} = {M'}_{\mkern-10mu X}^{\rm gp}$.
\end{enumerate}

\begin{proof}
Since $M_Y \to \mathcal O_Y$ is a logarithmic structure, the map $\varepsilon^{-1}(\mathcal O_Y^\ast) \to \mathcal O_Y^\ast$ is an isomorphism.  Since pushforward is left exact, this remains true after pushforward.  This proves~\ref{item:7}.

Since each $\tau : Y \to X$ is proper, $\tau_\ast \mathcal O_Y$ is coherent, hence a finite $\mathcal O_X$-algebra.  This proves~\ref{item:8}.

	All logarithmic modifications of $X$ are \'etale-locally pulled back from toric modifications of toric varieties and \ref{item:9} is a local property in the \'etale topology of $X$.  We therefore assume that $\tau$ is the strict base change of a toric modification $\tilde\tau : \tilde Y \to \tilde X$ is a toric modification.  We have $\tilde\tau_\ast(\mathcal O_{\tilde Y}) = \mathcal O_{\tilde X}$ (since $\tilde\tau$ is proper and birational and $\tilde X$ is normal, for example) so $\tilde\tau$ has geometrically connected fibers \cite[Th\'eor\`eme~(4.3.1) and Remarque~(4.3.4)]{ega3-1}, and therefore so does $\tau$.  This proves~\ref{item:9}.

If $x$ is a point of $X$ and $x'$ is its unique lift to $X'$ then the residue field $k(x')$ embeds in the residue field $K$ at the generic point of any irreducible component of $\tau^{-1} x$.  But the irreducible components are all rational, so $K$ is isomorphic to the field of rational functions (in some number of variables) over $k(x)$.  On the other hand, $k(x')$ is also an algebraic extension of $k(x)$, and the algebraic closure of $k(x)$ in the field of rational functions is $k(x)$.  Hence $k(x) = k(x')$.

We deduce~\ref{item:10} immediately from~\ref{item:2},~\ref{item:3}, and the topological invariance of the \'etale site~\cite[\href{https://stacks.math.columbia.edu/tag/04DZ}{Tag 04DZ}]{stacks-project}.

We prove~\ref{item:11}.  Apply the snake lemma to diagram~\eqref{eqn:a2} to get an exact sequence~\eqref{eqn:a3}:
\begin{gather} 
\label{eqn:a2} \vcenter{ \xymatrix{
0 \ar[r] & \mathcal O_X^\ast \ar[r] \ar[d] & M_X^{\rm gp} \ar[r] \ar[d] & \bar M_X^{\rm gp} \ar[r] \ar[d] & 0 \\
0 \ar[r] & \tau_\ast \mathcal O_Y^\ast \ar[r] & \tau_\ast M_Y^{\rm gp} \ar[r] & \tau_\ast \bar M_Y^{\rm gp} \ar[r] & \mathrm R^1 \tau_\ast \mathcal O_Y^\ast
} } \\
\label{eqn:a3}
	\tau_\ast(M_Y^{\rm gp}) / M_X^{\rm gp} \to \tau_\ast(\bar M_Y^{\rm gp}) / \bar M_X^{\rm gp} \xrightarrow{\phi} \mathrm R^1 \tau_\ast \mathcal O_Y^\ast
\end{gather}
	We saw in Lemma~\ref{lem:characteristic-injective} that $\phi$ is injective.  Indeed, $Y \to X$ is, \'etale-locally in $X$, the pullback of a toric modification $Y' \to X'$ and we can identify sections of $\tau_\ast( \bar M_Y^{\rm gp} ) / \bar M_X^{\rm gp}$ with the Picard group of $Y'$.  By Lemma~\ref{lem:characteristic-injective}, this injects into the Picard group of the fiber, and a fortiori into $\mathrm R^1 \tau_\ast \mathcal O_Y^\ast$.

	Thus the image of $\tau_\ast M_Y^{\rm gp}$ in $\tau_\ast \bar M_Y^{\rm gp}$ coincides with the image of $\bar M_X^{\rm gp}$.  That is, $\tau_\ast( M_Y^{\rm gp}  ) / \tau_\ast( \mathcal O_Y^\ast ) = \bar M_X^{\rm gp}$.  Furtheremore, the image of $\tau_\ast M_Y$ in $\bar M_X^{\rm gp}$ certainly contains $\bar M_X$ and is contained in $\tau_\ast(\bar M_Y) \cap \bar M_X^{\rm gp}$ (intersection taken in $\tau_\ast \bar M_Y^{\rm gp}$).  On the other hand, $\tau$ is a logarithmic modification, so it corresponds to a subdivision of a rational polyhedral cone $\sigma$, and sections of $\tau_\ast(\bar M_Y) \cap \bar M_X^{\rm gp}$ correspond to piecewise linear functions on that cone that are nonnegative on the subdivision.  But nonnegativity on a subdivision is equivalent to nonnegativity on the original cone, so we deduce that $\bar M_X = \tau_\ast(\bar M_Y^{\rm gp}) \cap \bar M_X$ and therefore that $\bar M'_X = \bar M_X$.  This proves~\ref{item:5}.

We note that the logarithmic \'etale site of $X$ depends only on the \'etale site and the \'etale sheaf $\bar M_X$.  Therefore~\ref{item:12} follows from~\ref{item:10} and~\ref{item:11}.

Finally,~\ref{item:13} follows immediately from~\ref{item:10} and the definition of $X'$.
\end{proof}

\paragraph{Lemma} \label{lem:log-mod-lim}
Let $X$ be a (fine and saturated) logarithmic scheme and let $\mathcal O'_X$ and $M'_X$ be the logarithmic \'etale sheafifications of $\mathcal O_X$ and $M_X$.  Let $X' = \Spec_X \mathcal O'_X$ with the logarithmic structure induced from $M'_X$.  Let $\mu : X' \to X$ be the projection.  Then we have the following properties:
\begin{enumerate}[label=(\roman{*})]
\item \label{item:1} $M'_X \to \mathcal O'_X$ is a logarithmic structure;
\item \label{item:2} the projection $\mu : X' \to X$ is integral (in the algebraic sense);
\item \label{item:3} $\mu$ is a universal bijection (in the category of schemes);
\item \label{item:16} the residue field extensions of $\mu : X' \to X$ are trivial;
\item \label{item:4} $\mu^\ast$ is an equivalence between the strict \'etale sites of $X$ and $X'$;
\item \label{item:5} $\bar M_X \to \bar M'_X$ is an isomorphism;
\item \label{item:6} $\mu^\ast$ is an equivalence between the logarithmic \'etale sites of $X$ and $X'$.
\item \label{item:14} $\mathrm R\mu_\ast \mathcal O_{X'} = \mathcal O'_X$ and $\mu_\ast M_{X'} = M'_X$ and $\mathrm R \mu_\ast M_{X'}^{\rm gp} = {M'}_{\mkern-10mu X}^{\rm gp}$.
\end{enumerate}

\begin{proof}
All assertions are local in the strict \'etale topology on $X$, so we assume that $X$ is quasicompact and quasiseparated and has a global chart.  

	By Corollary~\ref{cor:logification}, $M'_X = \varinjlim_{\tau} \tau_\ast M_Y$ and $\mathcal O'_X = \varinjlim_\tau \mathcal O_Y$, with the colimits taken over logarithmic modifications $\tau : Y \to X$. A filtered colimit of logarithmic structures is a logarithmic structure, so~\ref{item:7} of Lemma~\ref{lem:log-mod-pf} implies~\ref{item:1}.  A filtered colimit of finite algebras is integral, so~\ref{item:8} of Lemma~\ref{lem:log-mod-pf} implies~\ref{item:2}.  Formation of the characteristic monoid commutes with filtered colimits of logarithmic structures, so~\ref{item:11} of Lemma~\ref{lem:log-mod-pf} implies~\ref{item:5}.

If $\tau : Y \to X$ is a logarithmic modification and $x$ is a geometric point of $X$ then $k(x) \otimes \tau_\ast \mathcal O_Y$ is a $k(x)$-algebra with exactly one prime ideal, by \ref{item:9} of Lemma~\ref{lem:log-mod-pf}.  Equivalently, every non-unit is nilpotent.  This property is preserved by filtered colimits, so we have~\ref{item:3}.  A filtered colimit of isomorphisms of fields is an isomorphism, so~\ref{item:15} of Lemma~\ref{lem:log-mod-pf} implies~\ref{item:16}.

The topological invariance of the \'etale site gives~\ref{item:4} \cite[\href{https://stacks.math.columbia.edu/tag/04DZ}{Tag 04DZ}]{stacks-project}.

	Finally~\ref{item:6} follows from~\ref{item:4} and~\ref{item:5} and in Lemma~\ref{lem:log-mod-pf}; \ref{item:14} follows from~\ref{item:4} and the definition of $X'$.
\end{proof}

\paragraph{Corollary} \label{cor: logstrdescent=strdescent}
Let $X$ be a fine and saturated logarithmic scheme whose structure sheaf satisfies logarithmic \'etale descent. Then $M_X$ and $M_X^{\rm gp}$ also satisfy logarithmic \'etale descent.

\begin{proof}
	Let $\mu: X' \to X$ be the projection of Lemma~\ref{lem:log-mod-lim}.  By hypothesis, this is an isomorphism of underlying schemes.  It is also strict by \ref{lem:log-mod-lim}~\ref{item:5}, hence an isomorphism of logarithmic schemes.
\end{proof}

\section{Examples and counterexamples}
\label{sec:examples}

\paragraph{Lemma} \label{lem:toric-pushforward}
	Let $X$ be a strict closed subscheme of a toric variety defined by a monomial toric ideal.  Suppose that $\tau : Y \to X$ is a logarithmic modification.  Then $\mathcal O_X \to \tau_\ast \mathcal O_Y$ is surjective.

\begin{proof}
	By assumption, there is a toric modification $\tau' : Y' \to X'$ of which $\tau : Y \to X$ is the base change.  We have $\tau'_\ast \mathcal O_{Y'} = \mathcal O_{X'}$.  Let $T$ be the common torus of $X'$ and $Y'$.  Then the action of $T$ equips $\mathcal O_{X'}$ and $\mathcal O_{Y'}$ with a grading.  Since $X$ is defined in $X'$ by a monomial toric ideal, $Y$ is also defined inside $Y'$ by a monomial toric ideal.  A monomial toric ideal is a sum of graded pieces of $\mathcal O_{Y'}$ so $\mathcal O_Y$ is a direct summand of $\mathcal O_{Y'}$.  In particular, $\mathcal O_{X'} \simeq \tau'_\ast \mathcal O_{Y'} \to \tau_\ast \mathcal O_Y$ is surjective.  But this factors through $\mathcal O_X$ so $\mathcal O_X \to \tau_\ast \mathcal O_Y$ must also be surjective.
\end{proof}

In the following corollary, we say that a logarithmic scheme is \emph{idealized logarithmically flat} if, \'etale-locally in $X$, there is a toric variety $V$ with dense torus $T$, an equivariant subscheme $Z \subset V$, and a flat morphism $X \to [Z/T]$.  Equivalently, $X$ has a smooth cover $X' \to X$ (namely, $X' = X \mathop\times_{[V/T]} V$) and there is a flat morphism $X' \to Z$ where $Z$ is as above.

\paragraph{Corollary} \label{cor:1}
Suppose that $X$ is a fine and saturated logarithmic scheme that is idealized logarithmically flat.  Let $\tau : Y \to X$ be a logarithmic modification.  Then $\mathcal O_X \to \tau_\ast \mathcal O_Y$ is surjective.

\begin{proof}
	The hypotheses are preserved if $X$ is replaced by an fppf cover, and the conclusion is fppf-local on $X$, so we may assume that there is a flat morphism $X \to Z$ where $Z$ is an equivariant closed subscheme of a toric variety $V$.  Replacing $V$ by an equivariant open subvariety, we can assume that $V$ is affine and that the image of $X \to V$ meets the unique closed orbit.  Also replacing $X$ by an \'etale cover, we assume as well that the modification $\tau : Y \to X$ is pulled back from a toric modification $\tau' : W \to V$.  Then $\mathcal O_V \to \tau'_\ast \mathcal O_W$ is surjective by Lemma~\ref{lem:toric-pushforward}.  By flat base change, it follows that $\mathcal O_X \to \tau_\ast \mathcal O_Y$ is surjective as well.
\end{proof}

\paragraph{Corollary} \label{cor:2}
If $X$ is idealized logarithmically flat and reduced then $\mathcal O_X \to \tau_\ast \mathcal O_Y$ is an isomorphism.

\begin{proof}
	By Corollary~\ref{cor:2}, $\mathcal O_X \to \tau_\ast \mathcal O_Y$ is surjective.  On the other hand, Lemma~\ref{lem:log-mod-pf} implies that the kernel of $\mathcal O_X \to \tau_\ast \mathcal O_Y$ is a nil ideal.  Since $\mathcal O_X$ is reduced, the kernel must be zero.
\end{proof}

\paragraph{Corollary} \label{cor:field}
Suppose that $X$ is a fine and saturated logarithmic scheme whose underlying scheme is the spectrum of a field.  Let $\tau : Y \to X$ be a logarithmic modification.  Then $\tau_\ast \mathcal O_Y = \mathcal O_X$.

\begin{proof}
	With any logarithmic structure, the spectrum of a field is idealized logarithmically flat, and a field is certainly reduced, so the conclusion follows from Corollary~\ref{cor:2}.
\end{proof}

\paragraph{Example} \label{ex:nakayama-1}
If $\tau : Y \to X$ is a logarithmic modification, $\mathcal O_X \to \tau_\ast \mathcal O_Y$ can fail to be injective.  The following example is from \cite[p.\ 671, Remark]{Nakayama17}: 

Let $\mathbf A^2$ have its toric logarithmic structure and coordinates $x$ and $y$.  Let $X$ be the vanishing locus of $x^2$ and $y^2$ in $\mathbf A^2$, with logarithmic structure induced by restriction.  The element $xy$ is nonzero on $X$ but vanishes on the logarithmic blowup of $X$ at the toric ideal $(x,y)$.

\paragraph{Example}
We work over a field $k$. If the structure sheaf of a fine and saturated logarithmic scheme $X$ satisfies logarithmic \'etale descent and $Y \to X$ is a logarithmic modification, the structure sheaf of $Y$ need not necessarily satisfy logarithmic \'etale descent.  For an example, first let $\sigma$ be the rational polyhedral cone dual to the monoid $\mathbf N \alpha + \mathbf N \beta$ freely generated by $\alpha$ and $\beta$.  Let $\Delta$ be the subdivision of $\sigma$ along the ray $2 \alpha = 3 \beta$.  The associated morphism of toric varieties $\tilde Y \to \tilde X$ (equipped with their toric logarithmic structures) is a logarithmic modification.  Let $\tau \subset \Delta$ be the ray where $2 \alpha = 3 \beta$.  Its dual monoid is generated by the lattice $\mathbf Z(2 \alpha - 3 \beta)$ and the element $\alpha - \beta$.  Let $U \subset \tilde Y$ be the open subset associated with $\tau \subset \Delta$.  Then $U = \Spec k[x^{\pm (2 \alpha - 3 \beta)}, x^{\alpha - \beta}]$.  Since $x^\alpha = x^{3\beta - 2\alpha} (x^{\alpha-\beta})^3$ and $x^\beta = x^{3\beta - 2\alpha}(x^{\alpha-\beta})^2$, the fiber of $U$ over the origin $Z$, given by $x^\alpha = x^\beta = 0$, in $\tilde X$ is isomorphic to $\mathbf G_m \times \Spec k[\epsilon]/(\epsilon^2)$.  Its logarithmic structure is generated by $x^{\alpha-\beta}$, which corresponds to $\epsilon$.  

We can now consider $\tilde Y \times \tilde Y \to \tilde X \times \tilde X$.  The origin $X = Z \times Z \subset \tilde X \times \tilde X$ satisfies logarithmic \'etale descent by Corollary~\ref{cor:field} since it is the spectrum of a field.  On the other hand, the fiber $Y$ of $\tilde Y \times \tilde Y$ over $X$ contains the fiber of $U \times U$ over $Z$ as an open subset, and this open subset is isomorphic to $\mathbf Z[\epsilon_1, \epsilon_2] / (\epsilon_1^2, \epsilon_2^2)$ with its logarithmic structure generated by $\epsilon_1$ and $\epsilon_2$.  In the characteristic monoid, $\epsilon_1$ and $\epsilon_2$ correspond to $\alpha_1 - \beta_1$ and to $\alpha_2 - \beta_2$, respectively, where $\alpha_i$ and $\beta_i$ are pulled back from the coordinates on the two projections to $\tilde X$.  Thus Example~\eqref{ex:nakayama-1} appears on an open subscheme of a logarithmic blowup of $X$.  Moreover, the logarithmic blowup of $Y$ along the fractional ideal $(\alpha_1 - \beta_1, \alpha_2 - \beta_2)$ restricts to the logarithmic blowup of $Y \cap (U \times U)$ at $(\epsilon_1, \epsilon_2)$.  It follows that the structure sheaf of the logarithmic modification $Y$ of $X$ does not satisfy logarithmic \'etale descent.

\paragraph{Example}
If the structure sheaf of $X$ satisfies logarithmic \'etale descent and $X' \to X$ is a Kummer logarithmic \'etale morphism, the structure sheaf of $X'$ need not satisfy logarithmic \'etale descent.  For example, if $X$ is the spectrum of a field $k$ with a global chart by $\mathbf N^2$ then $X$ satisfies logarithmic \'etale descent by Corollary~\ref{cor:field}.  There is a Kummer logarithmically flat morphism $X' \to X$ where $X' = k[x,y] / (x^2, y^2)$, corresponding to the inclusion of monoids $\mathbf N^2  \subset \frac12 \mathbf N^2$.  If $k$ does not have characteristic~$2$ then this morphism is Kummer logarithmically \'etale.  We have seen in Example~\ref{ex:nakayama-1} that $X'$ does not satisfy logarithmic \'etale descent.

\paragraph{Example}
	If $\tau : Y \to X$ is a logarithmic modification, $\mathcal O_X \to \tau_\ast \mathcal O_Y$ can fail to be surjective.  The following example was communicated to us by C.\ Nakayama.  

	Let $X$ be defined by the ideal $(x^3 z - xy^2 z, xy^3 - xyz^2, yz^3 - x^2 yz)$ in $\Spec \mathbf Z[x,y,z]$.  Give $X$ the logarithmic structure freely generated by $x$, $y$, and $z$.  Let $Y$ be the logarithmic blowup of $X$ at the ideal generated by $(x,y,z)$.  Let $U \subset Y$ be the chart pulled back from $\Spec \mathbf Z[x,y/x,z/x]$, let $V \subset Y$ be the chart pulled back from $\Spec \mathbf Z[x/y,y,z/y]$, and let $W \subset Y$ be the chart pulled back from $\Spec \mathbf Z[x/z,y/z,z]$.  There is an element of $\Gamma(Y, \mathcal O_Y)$ defined locally by $x^4$ on $U$, by $y^4$ on $V$, and by $z^4$ on $W$.  Indeed, on $U \cap V$ we have $x^4 - x^2 y^2 = (x/y)^3 (x y^3 - x y z^2) +  (x/y) (z/y) (x^3 z - xy^2 z)$ so $x^4 \equiv x^2 y^2$ and, by symmetry, $x^2 y^2 \equiv y^4$.  Symmetric considerations show that $x^4 \equiv z^4$ on $U \cap w$ and $y^4 \equiv z^4$ on $V \cap W$.  

	On the other hand, suppose for the sake of contradiction that the element described above were the image of some putative $g \in \mathbf Z[x,y,z]$.  Two monomials of total degree~$4$ in $\mathbf Z[x,y,z]$ agree on $U \cap V \cap W$ if and only if the exponents of $x$, $y$, and $z$ in each have the same parity.  Therefore, $g$ must have an expression of the following form:
	\begin{equation*}
		g = a x^4 + b x^2 y^2 + c x^2 z^2 + d y^4 + e y^2 z^2 + f z^4
	\end{equation*}
	On the other hand, the elements of $\mathcal O_U$ with total degree $4$ and even exponents are freely generated by $x^4$, $x^2 y^2$, and $x^2 z^2 \equiv y^4 \equiv y^2 z^2 \equiv z^4$.  Therefore we must have $a = 1$, $b = 0$, and $c + d + e + f = 0$.  This equation and the symmetrical equations derived from restriction to $V$ and to $W$ are inconsistent, contradicting the assumed existence of $g$.

\paragraph{Lemma} \label{lem:alt}
If $S$ is a fine and saturated logarithmic scheme whose structure sheaf satisfies logarithmic \'etale descent and $\rho : T \to S$ is a logarithmic alteration (meaning a composition of logarithmic modifications and root stacks) then the natural map $\mathcal O_S \to \rho_\ast \mathcal O_{T}$ is an isomorphism.

\begin{proof}
	There is a logarithmic modification $\tau : U \to S$ and a root stack $\sigma : T \to U$ such that $\rho = \tau \sigma$.  We have $\sigma_\ast \mathcal O_T = \mathcal O_U$ by \cite[Theorem~3.1]{Kato21} and we have $\tau_\ast \mathcal O_U = \mathcal O_S$ by hypothesis.
\end{proof}

\paragraph{Proposition} \label{prop:flat-cover}
Let $\pi : X \to S$ be a morphism of fine and saturated logarithmic schemes. If the structure sheaf of $S$ satisfies logarithmic \'etale descent and $X$ is integral (in the logarithmic sense) and logarithmically flat over $S$ then the structure sheaf of $X$ satisfies logarithmic \'etale descent as well.

\begin{proof}
	By Lemma~\ref{lem:log-etale-sheaf}, it is sufficient to show that whenever $\tau : Y \to X$ is a logarithmic modification, there is a logarithmic alteration $\sigma : Y' \to X$ such that the natural map $\mathcal O_X \to \sigma_\ast \mathcal O_{Y'}$ is an isomorphism.  

	Note that $X$ is flat over $S$ because it is integral and logarithmically flat over $S$. Since $X$ is integral and saturated over $S$, the fine and saturated base change of any logarithmic alteration of $S$ coincides with the schematic base change.  Therefore it follows by flat base change that $\mathcal O_X \to \tau_\ast \mathcal O_{X'}$ is an isomorphism whenever $X'$ is the base change of a logarithmic alteration of $S'$.

	Consider a general logarithmic modification $\tau : Y \to X$.  By weak semistable reduction (~\cite[Theorem~1.0.1]{Mss}, \cite[Theorem 0.3]{AK}, \cite[Theorem 4.5]{ALT}), there are logarithmic alterations $\rho : S' \to S$ and $Y' \to Y$ fitting into a commutative square
	\begin{equation*} \begin{tikzcd}
		Y' \ar[r] \ar[d] & Y \ar[d] \\
		S' \ar[r,"\rho"] & S
	\end{tikzcd} \end{equation*}
	such that $Y'$ is integral and saturated over $S'$.  Let $X'$ be the base change of $X$ to $S'$ (the schematic base change is automatically integral and saturated since $X$ is integral and saturated over $S$).  Then $\tau' : Y' \to X'$ is a logarithmic modification of integral, logarithmically flat logarithmic schemes over $S'$, so by \cite[Proposition~4.4.13.1]{logpic}, we have $\mathrm R \tau'_\ast \mathcal O_{Y'} = \mathcal O_{X'}$.  Writing $\rho$ for the projection $X' \to X$ we have $\rho_\ast \mathcal O_{X'} = \mathcal O_X$ by Lemma~\ref{lem:alt}.  Therefore if $\sigma = \rho \tau' : Y' \to X$, we have $\mathrm R \sigma_\ast \mathcal O_{Y'}$, as required.
\end{proof}

\paragraph{Proposition} \label{prop:val-base}
If $X$ is a fine and saturated logarithmic scheme that is logarithmically flat over a fine and saturated logarithmic scheme whose characteristic monoid has rank~$1$ then the struture sheaf of $X$ satisfies logarithmic \'etale descent.

\begin{proof}
	Suppose that $\tau : Y \to X$ is a logarithmic modification.  Then $X$ and $Y$ are both integral over $S$ because all morphisms to a valuative logarithmic scheme are integral.  Therefore by \cite[Proposition~4.4.13.1]{logpic}, $\mathcal O_X \to \tau_\ast \mathcal O_Y$ is an isomorphism.  We conclude that the sructure sheaf of $X$ satisfies logarithmic \'etale descent by Lemma~\ref{lem:log-etale-sheaf}.
\end{proof}

The following proposition was suggested by C.\ Nakayama:

\paragraph{Proposition} \label{prop:log-reg}
If $X$ is a fine and saturated, logarithmically regular logarithmic scheme its structure sheaf satisfies logarithmic \'etale descent.

\begin{proof}
	Let $\tau : Y \to X$ be a logarithmic modification.  Then $Y$ is also logarithmically regular by \cite[Theorem~(8.2)]{Kato94}.  Therefore both $X$ and $Y$ are normal by \cite[Theorem~(4.1)]{Kato94}.  Moreover, the loci in $X$ and $Y$ where the logarithmic structures are trivial are dense in each, so $\tau : Y \to X$ is birational.  Thus $\tau$ is a proper, birational morphism between normal varieties, so $\mathcal O_X \to \tau_\ast \mathcal O_Y$ is an isomorphism.
\end{proof}

\section{Logarithmic \'etale descent}

First we consider Kummer logarithmic flat descent, following K.\ Kato \cite{Kato21}.

\paragraph{Lemma} \label{lem:root}
Let $\tau : Y \to X$ be a root stack.  Then $M_X \to \tau_\ast M_Y$ and $M_X^{\rm gp} \to \tau_\ast M_Y^{\rm gp}$ are isomorphisms, and $\mathrm R^1 \tau_\ast \bar M_Y^{\rm gp} = \mathrm R^1 \tau_\ast M_Y^{\rm gp} = 0$.

\begin{proof}
The maps $M_X \to \tau_\ast M_Y$ and $M_X^{\rm gp} \to \tau_\ast M_Y^{\rm gp}$ are isomorphisms by \cite[Theorem~3.1 and Theorem~3.2]{Kato21}.  The vanishing of $\mathrm R^1 \tau_\ast M_Y^{\rm gp}$ is essentially \cite[Corollary~5.1]{Kato21}, although the form of our statement is slightly more specific.  We give the details below, which will also imply the first two assertions.

Consider the following exact sequence:
\begin{equation*}
\tau_\ast(\bar M_Y^{\rm gp}) / \bar M_X^{\rm gp} \to \mathrm R^1 \tau_\ast \mathcal O_Y^\ast \to \mathrm R^1 \tau_\ast M_Y^{\rm gp} \to \mathrm R^1 \tau_\ast \bar M_Y^{\rm gp}
\end{equation*}
We argue that $\mathrm R^1 \tau_\ast \bar M_Y^{\rm gp} = 0$.  \'Etale-locally in $X$, it is possible to present $Y$ as the quotient of a finite $X$-scheme $Z$ by a finite diagonalizable group $G$.  Therefore $H^1(Y, \bar M_Y^{\rm gp})$ may be identified with the crossed homomorphisms $G \to H$, where $H(W) = \Gamma(W \mathop\times_X Z, \bar M_Z^{\rm gp})$ for flat $X$-schemes $W$, modulo principal crossed homomorphisms.  In this case, $\bar M_Y^{\rm gp}$ is pulled back from $X$, so the $G$-action is trivial, and crossed homomorphisms are the same as homomorphisms.  Finally $G$ is torsion and $\bar M_Z^{\rm gp}$ is torsion free and \'etale, so we conclude that $H^1(Y, \bar M_Y^{\rm gp}) = 0$.

It remains to show that $\tau_\ast(\bar M_Y^{\rm gp}) / \bar M_X^{\rm gp} \to \mathrm R^1 \tau_\ast \mathcal O_Y^\ast$ is an isomorphism.  We reduce this verifiation to the case where $X$ is the spectrum of a field.  First of all, by noetherian approximation, we can assume that $X$ is the spectrum of a noetherian ring.  The assertion is \'etale-local, so we can also assume $X$ is the spectrum of a henselian local ring.  Since giving an invertible sheaf, or an isomorphism between invertible sheaves, on $Y$ is a problem that is locally of finite presentation relative to $X$, we can apply Artin's approximation theorem to reduce to the case where $X$ is the spectrum of a complete noetherian local ring.  Then by Grothendieck's existence theorem it is sufficient to assume that $X$ is an infinitesimal extension of its closed point, $x$.  Any infinitesimal extension is an iteration of first-order extensions in which the ideal is isomorphic to $k(x)$.  The deformations and obstructions of invertible sheaves along such an extension are controlled by $H^1(Y, \mathcal O_Y)$ and $H^2(Y, \mathcal O_Y)$.  These groups may be interpreted as $G$-equivariant cohomology of $G$-linearized $\mathcal O_Z$-modules.  Since $G$ is diagonalizable, these groups vanish \cite[Th\'eor\`eme~5.3.3]{sga3-I}.  Thus if $X$ is an infinitesimal extension of a closed point $x$, then we have $H^1(Y, \mathcal O_Y^\ast) = H^1(\tau^{-1} x, \mathcal O_{\tau^{-1} x}^\ast)$.  It therefore suffices to prove the lemma in the case where $X = x$.

	Assume that $X$ is the spectrum of a separably closed field.  In this case, $Y$ has an explicit global description. The root stack corresponds to an extension of lattices $\bar M_X^{\rm gp} \rightarrow \bar N^{\rm gp}$ with finite cokernel. Let $\bar N$ be the saturation of $\bar M_X$ in $\bar N^{\rm gp}$. Then $Y$ is by definition the stack quotient of $Z = \Spec k[\bar N] \times_{\Spec k[\bar M_X]} \Spec k$ by the Cartier dual $G$ of $\bar N^{\rm gp}/ \bar M_X^{\rm gp}$. The logarithmic structure $M_Y$ on $Y$ is defined by descent. In particular, as $G$ acts trivially on $\bar N$, we have $\bar N = \tau_* \bar M_Y$. Note that $Z$ is a finite, connected (though possibly nonreduced) scheme over $X$, so it has trivial Picard group. Thus the Picard group of $Y$ consists of the $G$-linearizations of the trivial line bundle $\mathcal{O}_Z$, that is to say, the group of characters $\Hom(G, \mathbf{G}_{m,Z})$, where $\mathbf{G}_{m,Z}$ is the functor on $X$-schemes with values $\mathbf{G}_{m}(S \mathop\times_X Z) = \Gamma(S \mathop\times_X Z, \mathcal{O}_{S \times_X Z}^\ast)$.

As $Z$ is the spectrum of an artinian local $k$-algebra with residue field $k$, the group $\mathbf G_{m,Z}$ is a split extension of $\mathbf G_m$ by an iteratee extension of additive groups.  We can write $\mathbf G_{m,Z}$ as a product
\begin{equation} \label{eqn:93}
\mathbf{G}_{m,Z} \cong \mathbf{G}_{m} \times U 
\end{equation}
	with $U$ a unipotent group.  As $G$ is the Cartier dual of $\tau_\ast (\bar M_Y^{\rm gp})/ \bar M_X^{\rm gp}$, it is of multiplicative type, so there are no nontrivial homomorphisms from $G$ to $U$ by \cite[Proposition~2.4]{sga3-XVII}.  We therefore conclude
\begin{equation}
\Pic(Y) = \Hom(G,\mathbf{G}_{m,Z}) = \Hom(G,\mathbf{G}_{m}) = \tau_\ast (\bar M_Y^{\rm gp}) / \bar M_X^{\rm gp}
\end{equation}
by Cartier duality again.
\end{proof}

\paragraph{Lemma} \label{lem:saturated-hyperplane-logpic}
Suppose that $\tau : Y \to X$ is a morphism of logarithmic schemes, \'etale-locally (in $X$) pulled back from $\mathbf P^1 \to \logGm$ along a saturated morphism $X \to \logGm$.  Then every $M_Y^{\rm gp}$-torsor on $Y$ descends uniquely to a $\tau_\ast M_Y^{\rm gp}$-torsor on $X$.

\begin{proof}
	See \cite[Proposition~4.4.7]{logpic}.
\end{proof}

%

\paragraph{Lemma} \label{lem:hyperplane-subdivision}
If $\tau : Y \to X$ is a subdivision associated with $\alpha \in \bar M_X^{\rm gp}$ then $\mathrm R^1 \tau_\ast M_Y^{\rm gp} = 0$.

\begin{proof}
	By Lemma~\ref{lem:root-sat}, there is a root $\rho : X' \to X$ such that $n^{-1} \alpha$ is saturated.  Let $Y'$ be the base change, which is a root of $Y$.  Then $M_Y^{\rm gp} = \rho_\ast M_{Y'}^{\rm gp}$ and $\mathrm R^1 \rho_\ast M_{Y'}^{\rm gp} = 0$ by Lemma~\ref{lem:root}.  Therefore $\mathrm R^1 \tau_\ast M_Y^{\rm gp} = \mathrm R^1 (\tau \rho)_\ast M_{Y'}^{\rm gp}$.  On the other hand, we have $\mathrm R^1 \tau_\ast M_{Y'}^{\rm gp} = 0$ by Lemma~\ref{lem:saturated-hyperplane-logpic}.  Therefore $\mathrm R^1 (\tau \rho)_\ast M_{Y'}^{\rm gp} = \mathrm R^1 \rho_\ast ( \tau_\ast M_{Y'}^{\rm gp} )$.

It remains to show that every $\tau_\ast M_{Y'}^{\rm gp}$-torsor on $X'$ is \'etale-locally trivial on $X$.  This is certainly an \'etale-local question on $X$, so we assume that $X$ has a global chart by a monoid $P$ and $Y$ is the root associated with a Kummer extension $P \subset Q$.  Then we may present $X'$ as $[Z/G]$ where $G$ is the group that is Cartier dual to $Q^{\rm gp} / P^{\rm gp}$ and $\gamma : Z \to X$ is finite.  

We wish to show that any $\tau_\ast M_{Y'}^{\rm gp}$-torsor on $X'$ can be trivialized \'etale-locally in $X$.  By Lemma~\ref{lem:log-mod-pf}~\ref{item:11}, we have an exact sequence:
\begin{equation*}
0 \to \tau_\ast \mathcal O_{Y'}^\ast \to \tau_\ast M_{Y'}^{\rm gp} \to \bar M_{X'}^{\rm gp} \to 0
\end{equation*}
	We know that $\mathrm R^1 \rho_\ast \bar M_{X'}^{\rm gp} = 0$ from Lemma~\ref{lem:root}, so we will need to show that $\rho_\ast \bar M_{X'}^{\rm gp} \to \mathrm R^1 \rho_\ast( \tau_\ast \mathcal O_{Y'}^{\rm gp})$ is surjective.  We will actually prove the more precise statement that $\rho_\ast(\bar M_{X'}^{\rm gp}) / \bar M_X^{\rm gp} \to \mathrm R^1 \rho_\ast( \tau_\ast \mathcal O_{Y'}^\ast )$ is an isomorphism.  By noetherian approximation, it will suffice to prove this when $X$ is noetherian, and since it is \'etale-local, we can also replace $X$ with its henselization at a closed point.  By Artin's approximation theorem, we can replace $X$ by its completion, and by Grothendieck's existence theorem, we can assume further that $X$ is the spectrum of an artinian local ring.

	We may present $X'$ as $[Z/G]$ where $\gamma : Z \to X$ is finite and $G$ is a diagonalizable group.  We have a spectral sequence $H^p(G, R^q \gamma_\ast F \big|_Z) \Rightarrow \mathrm R^{p+q} \rho_\ast F$ for an \'etale sheaf $F$ on $X'$.  Since $Z$ is finite over $X$, the groups $\mathrm R^q \gamma_\ast F$ vanish for $q > 0$.  Therefore $\mathrm R^p \rho_\ast F = H^p(G, \gamma_\ast F)$ for all $p$.  In particular, $\mathrm R^1 \rho_\ast ( \tau_\ast \mathcal O_{Y'}^\ast )$ may be identified $H^1(G, \gamma_\ast \tau_\ast \mathcal O_{Y\times_X Z} )$.

	We assume that $X$ is artinian with closed point $x$.  Then $Z$ is also artinian, and so is $A = \mathcal O_Z \mathop\otimes_{\mathcal O_{X'}} \tau_\ast \mathcal O_{Y'} = \tau_\ast \mathcal O_{Z \mathop\times_X Y}$.  Let $A_0$ be the quotient of $A$ by its radical ideal.  Then $A^\ast$ is an iterated extension of $A_0^\ast$ by quasicoherent $\mathcal O_Z$-modules.  By \cite[Th\'eor\`eme~5.3.3]{sga3-I}, we have $H^p(G,J) = 0$ for all $p > 0$ and all $G$-linearized quasicoherent $\mathcal O_Z$-modules $J$, so by the long exact sequence in cohomology, $H^1(G, A^\ast) = H^1(G, A_0^\ast)$.

	Finally, we observe that $A_0$ is isomorphic to the quotient of $\mathcal O_Z$ by its nilradical.  In other words, its spectrum is the reduced structure on $\gamma^{-1} x$.  Then $H^1(G, A_0^\ast)$ is the same as $H^1(X'_0, \mathcal O_{X'_0}^\ast)$, where $X'_0$ is the reduced structure on $\tau^{-1} x$.  But $\tau^{-1}(x)^{\rm red} = \mathrm B G$ where $G$ is the group Cartier dual to $\rho_\ast (\bar M_{X'}^{\rm gp})_x / \bar M_{X,x}^{\rm gp}$.  In particular, its Picard group is isomorphic to the character group of $G$, which completes the proof.
\end{proof}

\paragraph{Theorem} \label{thm:log-mod-inv}
Suppose that $X$ is a fine and saturated logarithmic scheme.   Let $\tau : Y \to X$ be a logarithmic modification.  Write $M'_X$ and $M'_Y$ for the logarithmic \'etale sheafifications of $M_X$ and $M_Y$, respectively.  Then any ${M'}_{\mkern-10mu Y}^{\rm gp}$-torsor on $Y$ descends uniquely (up to unique isomorphism) to an ${M'}_{\mkern-10mu X}^{\rm gp}$-torsor on $X$.  In particular, $\tau_\ast {M'}_{\mkern-10mu Y}^{\rm gp} = {M'}_{\mkern-10mu X}^{\rm gp}$.
\begin{proof}
Let $\mu : X' \to X$ be the spectrum of the logarithmic \'etale sheafification of $\mathcal O_X$, equipped with the logarithmic structure induced from the logarithmic \'etale sheafification of $M_X$.  By Lemma~\ref{lem:log-mod-lim} the theorem will follow once it is proved with $X$ replaced with $X'$.  We therefore assume that $\mathcal O_X$ and $M_X$ are logarithmic \'etale sheaves.

The hypotheses are preserved by \'etale localization on $X$, and the conclusions are \'etale-local as well, so we assume that $X$ is quasicompact and quasiseparated and has a global chart.  Then $Y \to X$ is quasicompact and quasiseparated.  By Lemma~\ref{lem:sheafification}, we therefore have ${M'}_{\mkern-10mu Y}^{\rm gp} = \varinjlim \rho_\ast M_Z^{\rm gp}$ with the colimit taken over all logarithmic modifications $\rho : Z \to Y$.  We will argue that any ${M'}_{\mkern-10mu Y}^{\rm gp}$-torsor $L$ on $Y$ descends uniquely to an ${M'}_{\mkern-10mu X}^{\rm gp}$-torsor on $X$.

First we note that it is sufficient to assume that $L$ is induced from an $M_Y^{\rm gp}$-torsor on $Y$.  Indeed, suppose that the theorem is proved under this assumption.  If $L$ is any ${M'}_{\mkern-10mu Y}^{\rm gp}$-torsor on $Y$ then, because $Y$ is quasicompact and quasiseparated, and ${M'}_{\mkern-10mu Y}^{\rm gp} = \varinjlim \rho_\ast M_Z^{\rm gp}$, there is some logarithmic modification $\rho : Z \to Y$ and a $\rho_\ast M_Z^{\rm gp}$-torsor $L'$ that induces $L$.  Then $\rho^{-1} L' \otimes_{\rho^{-1} \rho_\ast M_Z^{\rm gp}} M_Z^{\rm gp}$ is an $M_Z^{\rm gp}$-torsor inducing the ${M'}_{\mkern-10mu Z}^{\rm gp}$-torsor $\rho^\ast L = \rho^{-1} L \otimes_{\rho^{-1} {M'}_{\mkern-10mu Y}^{\rm gp}} {M'}_{\mkern-10mu Z}^{\rm gp}$.  By our assumption, $\rho^\ast L$ descends uniquely to $X$ and to $Y$.  Let $K = \tau_\ast \rho_\ast \rho^\ast L$ be the unique descent to $X$.  Then $\tau^\ast K$ and $L$ both pull back to $\rho^\ast L$ on $Z$, so by the uniqueness part of our hypothesis, we must have a unique isomorphism $\tau^\ast K \simeq L$.  We therefore assume that $L$ is induced from an $M_Y^{\rm gp}$-torsor from now on.  

Next, we observe that there is an iterated subdivision $\rho : Z \to Y$ of $Y$ along hyperplanes such that $\tau \rho : Z \to X$ is also an iterated subdivision of $X$ by hyperplanes.  By the same argument as in the last paragraph, it will suffice to demonstrate the theorem under the additional assumption that $Y$ is an iterated subdivision of $X$ along hyperplanes.  We assume this as well from now on.

Suppose that $Y = Y_n \xrightarrow{\tau_n} Y_{n_1} \xrightarrow{\tau_{n-1}} \cdots \xrightarrow{\tau_1} Y_0 = X$ is a sequence of hyperplane subdivisions and $L$ is an $M_Y^{\rm gp}$-torsor.  For each $i$, let $Y'_i$ be the Stein factorization of $Y \to Y'_i$, equipped with the logarithmic structure pushed forward from $Y$.  We therefore obtain a sequence $Y = Y'_n \xrightarrow{\tau'_n} \cdots \xrightarrow{\tau'_1} Y'_0 = X$ (note that $Y'_0 = X$ because we assumed $\mathcal O_X$ is a logarithmic \'etale sheaf) with the property that ${\tau'_i}_\ast M_{Y_i}^{\rm gp} = M_{Y_{i-1}}^{\rm gp}$.  We now apply Lemma~\ref{lem:hyperplane-subdivision} inductively to descend an $M_{Y_i}^{\rm gp}$-torsor uniquely to a ${\tau'_i}_\ast M_{Y_i}^{\rm gp} = M_{Y_{i-1}}^{\rm gp}$-torsor and the proof is complete.  
\end{proof}

\paragraph{Corollary}
Suppose $X$ is a fine and saturated logarithmic scheme that is logarithmically flat over a fine and saturated logarithmic scheme whose characteristic monoid has rank~$\leq 1$. If $Y \to X$ is a logarithmic modification then 
\begin{equation*}
H^i(X,M_X^{\rm gp}) = H^i(Y,M_Y^{\rm gp})
\end{equation*}
for $i=0,1$.

\begin{proof}
	Since $X$ is logarithmically flat over a rank~$1$ base, so is $Y$, and thus the logarithmic structures of both $X$ and $Y$ satisfy logarithmic \'etale descent by Proposition~\ref{prop:val-base}. Thus $M_X' = M_X$, $M_Y'=M_Y$, and so the result follows from Theorem~\ref{thm:log-mod-inv}. 
\end{proof}

\section{The logarithmic Picard group}

Let $S$ be a fine and saturated logarithmic scheme and let $X$ be an integral, saturated, proper, vertical, logarithmically smooth curve over $S$.  In \cite{logpic}, the logarithmic Picard stack $\bLogPic(X/S)$ is defined to be the stack in the strict \'etale topology on the category $\bLogSch/S$ of fine and saturated logarithmic schemes over $S$ whose fiber over $S'$ is the category of $M_{X'}^{\rm gp}$-torsors of \emph{bounded monodromy} on $X'$, where $X'$ is the base change of $X$ to $S'$.  The sheaf (in the strict \'etale topology of $\bLogSch/S$) of isomorphism classes in $\bLogPic(X/S)$ is denoted $\LogPic(X/S)$.

The bounded monodromy condition is important for the existence of a logarithmically smooth cover of $\bLogPic(X/S)$, and of a logarithmically \'etale cover of $\LogPic(X/S)$, by a logarithmic scheme.  All we will need to know about bounded monodromy here is the following lemma, which can be deduced formally from \cite{logpic}; we therefore will not define bounded monodromy explicitly.

\paragraph{Lemma} \label{lem:bdd-mono}
Suppose that $X$ is an integral, saturated, proper, vertical, logarithmically smooth curve over $S$.  Let $L$ be a $M_X^{\rm gp}$-torsor on $X$.  Let $S' \to S$ be a logarithmic modification and let $X'$ be the base change of $X$ to $S'$.  Let $L'$ be the pullback of $L$ to $X'$.  If $L'$ has bounded monodromy then so does $L$.

\begin{proof}
	By \cite[Proposition~4.3.2]{logpic}, bounded monodromy can be verified at the valuative geometric points of $S$.  But every valuative geometric point of $S$ lifts to $S'$, since valuative logarithmic schemes have no nontrivial logarithmic modifications.
\end{proof}

\paragraph{Corollary} \label{cor:logpic}
Suppose that $X$ is a saturated, proper, vertical, logarithmically smooth curve over $S$.  Then $\bLogPic(X/S)$ forms a stack, and $\LogPic(X/S)$ forms a sheaf, in the logarithmic \'etale topology on the subcategory of $\bLogSch/S$ consisting of fine and saturated logarithmic schemes whose structure sheaves satisfy logarithmic \'etale descent.

\begin{proof}
	By Lemma~\ref{lem:log-etale-sheaf}, we must show that if $S$ is a logarithmic scheme whose structure sheaf satisfies logarithmic \'etale descent, $\tau : S' \to S$ is a logarithmic modification, and $\tau : X' \to X$ is the base change of $X$ along $\tau$, then the pullback morphisms $\bLogPic(X) \to \bLogPic(X')$ and $\LogPic(X) \to \LogPic(X')$ are equivalences.  

	The latter statement follows from the former: It is \'etale local in $S$, so we can assume that $X$ has a section through each of its connected components. This allows us to identify $\bLogPic(X/S) \simeq \LogPic(X/S) \times \mathrm B \logGm^r$, where $r$ is the (locally constant) number of connected components of the fibers of $X$ over $S$.  If $\bLogPic(X) \to \bLogPic(X')$ is known to be an equivalence then $\LogPic(X) \times \mathrm B \logGm^r(S) \to \LogPic(X') \times \mathrm B \logGm^r(S')$ is an equivalence.  On the other hand, $\mathrm B \logGm^r(S) \to \mathrm B \logGm^r(S')$ is an equivalence by Theorem~\ref{thm:log-mod-inv}.  Since $\mathrm B \logGm^r(S) \neq \varnothing$, we conclude that $\LogPic(X) \to \LogPic(X')$ is an isomorphism.

	We show now that $\bLogPic(X) \to \bLogPic(X')$ is an equivalence.  An object of $\bLogPic(X')$ is a $M_{X'}^{\rm gp}$-torsor $L'$ on $X'$ with bounded monodromy.  Since $X'$ is the base change of a logarithmic modification of $S$, it is a logarithmic modification of $X$.  By Proposition~\ref{prop:flat-cover}, the structure sheaf of $\mathcal O_X$ satisfies logarithmic \'etale descent.  Therefore $L'$ descends uniquely to a $M_X^{\rm gp}$-torsor $L$ on $X$ by Theorem~\ref{thm:log-mod-inv}.  It has bounded monodromy by Lemma~\ref{lem:bdd-mono}.
\end{proof}

\bibliographystyle{halpha}
\bibliography{refs,me,ega-sga}

\end{document}